\documentclass[12pt]{article}
\usepackage[cp1251]{inputenc}
\usepackage[english]{babel}
\usepackage{amsfonts}

\newtheorem{lemma}{Lemma}
\newtheorem{theorem}{Theorem}

\newtheorem{corollary}{Corollary}

\begin{document}

\title{Spectral data for Hamiltonian-minimal Lagrangian tori in ${\mathbb C}P^2$}
\author{A. E. Mironov
\thanks{
This work was supported by the RFBR (grant 06-01-00094a), by the
SB RAS (complex integration project No. 2.15), and by IHES,
France}}
\date{}

\maketitle
\section{Introduction}
In this work, we find spectral data that allow to find
Hamiltonian-minimal Lagrangian tori in ${\mathbb C}P^2$ in terms
of theta functions of spectral curves.

A Lagrangian submanifold in a K\"ahler manifold is called
Hamiltonian-minimal, if the variations of volume along the
Hamiltonian fields are equal to zero.

The simplest example of a Hamiltonian-minimal Lagrangian
submanifold ({\bf HML}-submanifold) is a Clifford torus \cite{O}
$$
 S^1(r_1)\times\dots\times S^1(r_n)\subset {\mathbb C}^n,
$$
where $S^1(r_i)\subset{\mathbb C}$ is a circle of radius $r_i$.
See \cite{M1} for more examples of closed {\bf HML}-submanifolds
in ${\mathbb C}^n$ and ${\mathbb C}P^n$.

Helein and Romon \cite{HR1} give a description of {\bf HML}-tori in ${\mathbb
C}^2$ using a Weierstrass representation of Lagrangian surfaces. McIntosh and Romon \cite{MR} find spectral data for such surfaces.

Until recently, only a few examples of Hamiltonian-minimal (but
not minimal) Lagrangian tori in ${\mathbb C}P^2$ were known (see
\cite{M1}, \cite{M2}, \cite{HM1}). The authors of \cite{HR2},
\cite{HM2} show that {\bf HML}-tori in ${\mathbb C}P^2$ are finite
type surfaces.

Below we show that {\bf HML}-tori in ${\mathbb C}P^2$ are described by the following system of equations (Lemma 3)
$$
 V_y+U_x=0,
$$
$$
2U_y-2V_x=(bv_x+av_y)e^v,
$$
$$
  \partial_x^2 v+\partial_y^2 v=4(U^2+V^2)e^{-2v}-4e^{v}-2(Ua-Vb)e^{-v},
$$
where $ds^2=2e^{v(x,y)}(dx^2+dy^2)$ is an induced metric on a
torus, $U(x,y)$, $V(x,y)$ are periodic real functions, $a,b$ are
some real constants. We do not know whether this system has a Lax
representation with a spectral parameter. This is the main
difficulty for constructing quasi-periodic solu\-tions in terms of
theta functions of spectral curves.

Minimal Lagrangian tori correspond to the condition $a=b=0$. In
this case, the indicated system of equations is reduced to the
Tzizeica equation
$$
\partial_x^2v+\partial_y^2v=4e^{-2v}-4e^v.
$$
The Tzizeica equation allows the Lax representation with a spectral parameter found in \cite{Mih}. Quasi-periodic solutions of this equation were found in \cite{Sh}. In \cite{Mc}, the authors show that there are smooth periodic solutions corresponding to smooth minimal Lagrangian tori.

We briefly present the basic elements of our construction.

Define a Lagrangian surface $\Sigma\subset{\mathbb C}P^2$ using a composition of the mappings
$$
 \varphi:{\mathbb R}^2\rightarrow S^5\subset {\mathbb C}^3
$$
and the Hopf projection
$$
 {\cal H}:S^5\rightarrow {\mathbb C}P^2,
$$
where $S^5$ is a unite sphere in ${\mathbb C}^3$. Suppose that the
induced metric on $\Sigma$ has a conformal form
$$
 ds^2=2 e^{v(x,y)}(dx^2+dy^2),
$$
where $x,y$ are coordinates on ${\mathbb R}^2$. As follows from the
Lagrangianity and the conformality of mapping ${\cal H}\circ\varphi$,
$$
 <\varphi,\varphi_x>=<\varphi,\varphi_y>=<\varphi_x,\varphi_y>=0,\
 |\varphi_x|^2=|\varphi_y|^2=2e^{v}, \eqno{(1)}
$$
where $<.,.>$ is a Hermitian product in ${\mathbb C}^3$ (see
\cite{M2}). We are looking for the components $\varphi^i, i=1,2,3$
of the vector function $\varphi$ in the form
$$
 \varphi^i=\alpha_i\psi(x,y,Q_i),
$$
where $\psi(x,y,Q)$ is a two-points Baker-Akhiezer function on an
auxiliar Rie\-mannian surface $\Gamma$ named {\it spectral curve},
$Q_i\in\Gamma$ are some points, $\alpha_i$ --- some constants. In
theorem 1, we indicate spectral data, i.e. some sufficient
constraints on the Baker-Akhiezer function $\psi$, such that the
equations
$$
 <\varphi,\varphi>=1,\ <\varphi,\varphi_x>=<\varphi,\varphi_y>=<\varphi_x,\varphi_y>=0
$$
hold.

Note that equation
$$
 <\varphi_x,\varphi_y>=
 \varphi^1_x\bar{\varphi}^1_y+
 \varphi^2_x\bar{\varphi}^2_y+
 \varphi^3_x\bar{\varphi}^3_y=0,
$$
is "similar" to equation
$$
 \frac{\partial x^1}{\partial u^1}\frac{\partial x^1}{\partial u^2}+
 \frac{\partial x^2}{\partial u^1}\frac{\partial x^2}{\partial u^2}+
 \frac{\partial x^3}{\partial u^1}\frac{\partial x^3}{\partial u^2}=0,
$$
that needs to be solved in order to construct orthogonal curvilinear systems of coordinates in
${\mathbb R^3}$.
Krichever \cite{K} indicates spectral data allowing to restore $n$-orthogonal curvilinear coordinates in ${\mathbb R^n}$. The idea of theorem 1 is close to \cite{K}.

In lemmas 5 and 6 we indicate constraints on spectral data that
let mapping ${\cal H}\circ\varphi$ be conformal, i.e. condition
$$
|\varphi_x|^2=|\varphi_y|^2=2e^{v}
$$
holds.

The key point of our construction is assertion (lemma 2) stating
that the components of the vector function $\varphi$ satisfy the
Schr\"odinger equation
$$
 \partial_x^2\varphi^j+\partial_y^2\varphi^j-
 i(\beta_x\partial_x\varphi^j+\beta_y\partial_y\varphi^j)+4e^v\varphi^j=0,
$$
where $\beta$ is a Lagrangian angle (for definition see below). The Lagrangian angle defines the mean curvature vector
$$
 H=J\nabla \beta,
$$
$\nabla\beta$ being the gradient of function $\beta$ in the induced metric on a torus, $J$ being the complex structure operator on ${\mathbb C}P^2$.

In case that surface $\Sigma$ is a {\bf HML}-torus, the Lagrangian
angle is a linear function, i.e. $\beta_x$ and $\beta_y$ are
constants (see \cite{M1}). If surface $\Sigma$ is a minimal
Lagrangian torus, the Lagrangian angle is constant, i.e. the
components $\varphi^j$ satisfy the potential Schr\"odinger
equation
$$
 \partial_x^2\varphi^j+\partial_y^2\varphi^j+4e^v\varphi^j=0.
$$
Thus, the Schr\"odinger equation allows to distinguish {\bf
HML}-tori and minimal Lagrangian tori amongst Lagrangian tori. In
the main theorem 2 we indicate spectral data corresponding to such
surfaces.

The method used to distinguish Lagrangian surfaces with a
conformal metric (lemma 5) was first applied by Novikov and
Veselov \cite{VN} for distingu\-ishing two-dimensional potential
Schr\"odinger operators that are finite-gap on one level of energy
from finite-gap Schr\"odinger operators with a magnetic field. All
in all, the spectral data found in theorem 2 are a generalization
of spectral data found in \cite{VN}.

The rest of this paper is organized as follows. In paragraph 2 we
find equations for {\bf HML}-tori in ${\mathbb C}P^2$. In
paragraph 3 we prove the main theorem. In paragraph 4 we give
examples of spectral curves that satisfy the conditi\-ons stated
in the main theorem.

Note that our solutions are, in the general case, quasi-periodic, and we do not ask whether there exist periodic solutions. This problem requires extra research which we are planning to conduct in the future.

This work is dedicated to Sergey Petrovich Novikov on the occasion of his 70th birthday.

\section{Equations of Hamiltonian-minimal Lagrangian tori in ${\mathbb C}P^2$}
As follows from (1), matrix
$$
 \tilde{\Phi}=
 (\varphi,\frac{1}{\sqrt{2}}e^{-\frac{v}{2}}\varphi_x,
 \frac{1}{\sqrt{2}}e^{-\frac{v}{2}}\varphi_y)^{\top}
$$
belongs to the group $U(3)$. A {\it Lagrangian angle} is a function defined from the equality
$$
 e^{i\beta (x,y)}={\rm det} \tilde{\Phi}.
$$
From the definition of a Lagrangian angle we get
$$
 \Phi=
  \left(
  \begin{array}{ccc}
    \varphi^1 & \varphi^2 &  \varphi^3\\
   \frac{1}{\sqrt{2}}e^{-\frac{v}{2}-i\frac{\beta}{2}}\varphi^1_x
   &  \frac{1}{\sqrt{2}}e^{-\frac{v}{2}-i\frac{\beta}{2}}\varphi^2_x &
    \frac{1}{\sqrt{2}}e^{-\frac{v}{2}-i\frac{\beta}{2}}\varphi^3_x \\
    \frac{1}{\sqrt{2}}e^{-\frac{v}{2}-i\frac{\beta}{2}}\varphi^1_y
    &  \frac{1}{\sqrt{2}}e^{-\frac{v}{2}-i\frac{\beta}{2}}\varphi^2_y
    &  \frac{1}{\sqrt{2}}e^{-\frac{v}{2}-i\frac{\beta}{2}}\varphi^3_y \\
  \end{array}\right)\in{\rm SU(3)},
$$
Matrix $\Phi$ satisfies equations
$$
 \Phi_x=A\Phi,\ \Phi_y=B\Phi,\eqno{(2)}
$$
where $A$ and $B$ have the form
$$
 A=
  \left(
  \begin{array}{ccc}
   0 & \sqrt{2}e^{\frac{v}{2}+i\frac{\beta}{2}} & 0\\
  -\sqrt{2}e^{\frac{v}{2}-i\frac{\beta}{2}} & if &
  -\frac{v_y}{2}+i(h+\frac{\beta_y}{2}) \\
  0 & \frac{v_y}{2}+i(h+\frac{\beta_y}{2}) & -if\\
  \end{array}\right)\in{\rm su(3)},
$$
$$
 B=
  \left(
  \begin{array}{ccc}
   0 & 0 & \sqrt{2}e^{\frac{v}{2}+i\frac{\beta}{2}} \\
  0 & ih & \frac{v_x}{2}+i(-f+\frac{\beta_x}{2})  \\
  -\sqrt{2}e^{\frac{v}{2}-i\frac{\beta}{2}} &
  -\frac{v_x}{2}+i(-f+\frac{\beta_x}{2}) & -ih\\
  \end{array}\right)\in{\rm su(3)},
$$
$f(x,y)$ and $h(x,y)$ being some real functions.
The zero-curvature equation
$$
 A_y-B_x+[A,B]=0
$$
 implies the following lemma (see \cite{M3})
\begin{lemma}
Lagrangian surfaces are defined by a system of equations
$$
2V_y+2U_x=(\beta_{xx}-\beta_{yy})e^v,
$$
$$
2U_y-2V_x=(\beta_yv_x+\beta_xv_y)e^v,
$$
$$
\Delta v=4(U^2+V^2)e^{-2v}-4e^{v}-2(U\beta_x-V\beta_y)e^{-v},
$$
where $U=fe^{v},V=he^{v}.$
\end{lemma}
From (2), get equalities
$$
 \varphi_{xx}^j=\frac{1}{2}(-4e^{v}\varphi^j+\varphi_x^j(2if+v_x+i\beta_x)+\varphi_y^j(2ih-v_y+i\beta_y)),
$$
$$
 \varphi_{yy}^j=\frac{1}{2}(-4e^{v}\varphi^j+\varphi_x^j(-2if-v_x+i\beta_x)+\varphi_y^j(-2ih+v_y+i\beta_y)),
$$
implying lemma (see \cite{M3}).
\begin{lemma}
The components $\varphi^j$ of the vector function $\varphi$
satisfy the Schr\"odin\-ger equation
$$
 \partial_x^2\varphi^j+\partial_y^2\varphi^j-i(\beta_x\partial_x\varphi^j+\beta_y\partial_y\varphi^j)+4e^v\varphi^j=0.
$$
\end{lemma}
Thus, there is a two-dimensional Schr\"odinger operator connected
in a natural way to every Lagrangian surface given by the mapping
$\varphi$. In of \cite{M3} (see also \cite{M4}) it is shown that
in case of minimal tori, the isospectral deformations of the
Schr\"odinger operator given by the Veselov-Novikov equation
correspond to integrable deformations of minimal tori.

If $\beta$ is a harmonic function, then mapping ${\cal H}\circ\varphi$
defines the {\bf HML}-surface \cite{M1}. In particular, if mapping ${\cal H}\circ\varphi$ is doubly periodic, then the Lagrangian angle is the linear function
$$
 \beta(x,y)= a x+ b y+ c,\
 a,b,c\in{\mathbb R}.
$$
From Lemma 1 follows
\begin{lemma}
 {\bf HML}-tori in ${\mathbb C}P^2$ are described by the following equations
 $$
   V_y+U_x=0,
 $$
 $$
  2U_y-2V_x=(bv_x+av_y)e^v,
 $$
 $$
  \Delta v=4(U^2+V^2)e^{-2v}-4e^{v}-2(Ua-Vb)e^{-v}.
 $$
\end{lemma}

Note that after replacing $V=-\frac{u_x}{2}, U=\frac{u_y}{2}$
the equations in lemma 3 reduce to a system on functions $u$ and $v$
$$
 \Delta u=e^v\nabla v A,
$$
$$
  \Delta v=e^{-2v}(\nabla u)^2-4e^{v}-e^{-v}\nabla u A,
$$
$A=(b,a)$, that passes the Painleve test for integrability such
that it could well be that there exists a Lax representation with
a spectral parameter for the equations of lemma 3. The fact that
there exist solutions in theta functions points into the same
direction.

\section{Main Theorem}
\subsection{Baker-Akhiezer function}

Let $\Gamma$ be a Riemannian surface of genus $g$.
Suppose that the following set of divisors be given on $\Gamma$
$$
 \gamma=\gamma_1+\dots+\gamma_{g+l},\
 r+r_1+\dots+r_l,
$$
a pair of marked points $P_1$, $P_2\in \Gamma$, and also the local
parameters $k_1^{-1}, k_2^{-1}$ in the neighborhood of points
$P_1$ and $P_2$. The {\it two-points Baker-Akhiezer function},
corresponding to the {\it spectral data}
$$
 \{\Gamma,P_1,P_2,k_1,k_2,\gamma,r+r_1+\dots+r_l\},
$$
 is called
a function $\psi(x,y,P), P\in\Gamma$, with the following
characteristics:

{\bf 1)} in the neighborhood of $P_1$ and $P_2$ function $\psi$
has essential singularities of the following form
$$
 \psi=e^{k_1x}\left(f_1(x,y)+\frac{g_1(x,y)}{k_1}+\frac{h_1(x,y)}{k_1^2}+\dots\right),
$$
$$
 \psi=e^{k_2y}\left(f_2(x,y)+\frac{g_2(x,y)}{k_2}+\frac{h_2(x,y)}{k_2^2}+\dots\right).
$$

{\bf 2)} on $\Gamma\backslash \{P_1,P_2\}$, function $\psi$
is meromorphic with simple poles on $\gamma$.

The space of such functions has dimension $l+1$, consequently,
the Baker-Akhiezer function is uniquely defined if we demand the following normalization conditions to be fulfilled

{\bf 3)}
$$
\psi(x,y,r)=d,\ \psi(x,y,r_i)=0, \ i=1,\dots,l,
$$
where $d$ is a non-zero constant.

Express the Baker-Akhiezer function explicitly by the theta
function of surface $\Gamma$.

On surface $\Gamma$, choose a base of cycles
$$
 a_1,\dots,a_g, \ b_1,\dots,b_g
$$
with the following intersection indices
$$
 a_i\circ a_j=b_i\circ b_j=0,\ a_i\circ b_j=\delta_{ij}.
$$
By $\omega_1,\dots,\omega_g$ we denote a base if holomorphic differentials on $\Gamma$, normalized by conditions
$$
 \int_{a_j}\omega_i=\delta_{ij}.
$$
Let the matrix of $b$-periods of differentials $\omega_j$ with components
$$
 B_{ij}=\int_{b_i}\omega_j
$$
be denoted by $B$. This matrix is symmetric and has a positively defined imaginary part.

The Riemannian theta function is defined by an absolutely
converging series
$$
 \theta(z)=\sum_{m\in{\mathbb Z}}e^{\pi i (Bm,m)+2\pi i(m,z)},\
 z=(z_1,\dots,z_g)\in{\mathbb C}^g.
$$
The theta function has the following characteristics
$$
 \theta(z+m)=\theta(z),\ m\in{\mathbb Z},
$$
$$
 \theta(z+Bm)={\rm exp}(-\pi i (Bm,m)-2\pi i (m,z))\theta(z),\ m\in{\mathbb
 Z}.
$$
Let $X$ denote a Jacoby variety of surface $\Gamma$
$$
 X={\mathbb C}^g/ \{{\mathbb Z}^g+B{\mathbb Z}^g\}.
$$
Let $A:\Gamma\rightarrow X$ be an Abelian mapping defined by formula
$$
 A(P)=\left(\int_{q_0}^P\omega_1,\dots,\int_{q_0}^P\omega_g\right),\
 P\in\Gamma,
$$
$q_0\in\Gamma$ being a fixed point.

For points $\gamma_1,\dots,\gamma_g$ in the general case and according to the Riemann theorem, the function
$$
 \theta(z+A(P)),
$$
where $z=K-A(\gamma_1)-\dots-A(\gamma_g)$, has exactly $g$ zeros on $\Gamma$ $\gamma_1,\dots,\gamma_g$.

Let $\Omega^1$ and $\Omega^2$ denote meromorphic differentials on $\Gamma$ with the only simple poles in points $P_1$ and
$P_2$, respectively and normalized by the conditions
$$
 \int_{a_j}\Omega^1=\int_{a_j}\Omega^2=0,\ j=1,\dots,g.
$$
Let $U$ and $V$ denote their vectors of $b$-periods
$$
 U=\left(\int_{b_1}\Omega^1,\dots,\int_{b_g}\Omega^1\right),\
 V=\left(\int_{b_1}\Omega^2,\dots,\int_{b_g}\Omega^2\right).
$$
To begin with, consider the case $l=0$.

Let $\tilde{\psi}$ denote function
$$
 \tilde{\psi}(x,y,P)=\frac{\theta(A(P)+xU+yV+z)}{\theta(A(P)+z)}
 {\rm exp}(2\pi ix\int_{q_0}^P\Omega^1+2\pi iy\int_{P_1}^P\Omega^2).\eqno{(3)}
$$
Then, the Baker-Akhiezer function sought for has the form
$$
 \psi(x,y,P)=f(x,y)\tilde{\psi}(x,y,P),
$$
where function $f$ is defined by the condition $\psi(x,y,r)=d$.

For $l>0$, the Baker-Akhiezer function can be found in the
foll\-ow\-ing form
$$
 \psi=f\tilde{\psi}+f_1\tilde{\psi}_1+\dots+f_l\tilde{\psi}_l,
$$
where $\tilde{\psi}_j$ is a function of form (3), constructed by the divisor $$
 \gamma_1+\dots+\gamma_{g-1}+r_j.
$$
Find functions $f,f_j$ from the normalization conditions 3).

\subsection{Lagrangian immersions }
Let $\varphi^1,\varphi^2,\varphi^3$ denote the following functions
$$
 \varphi^i=\alpha_i\psi(x,y,Q_i),
$$
where $Q_1,Q_2,Q_3\in\Gamma$ is an additional set of points, $\alpha_i$
are some constants.

In the following theorem, we state sufficient conditions for the
vector-function $\varphi=(\varphi^1,\varphi^2,\varphi^3)$ to
define a Lagrangian imbedding of the plane ${\mathbb R}^2$ in
${\mathbb C}P^2$.

Suppose that surface $\Gamma$ has a holomorphic involution
$\sigma$ and an antiholo\-morphic involution $\mu$
$$
 \sigma:\Gamma\rightarrow\Gamma,\ \mu:\Gamma\rightarrow\Gamma,
$$
for which points $P_1$, $P_2$ and $r$ are fixed, and
$$
 k_i(\sigma(P))=-k_i(P),\
 k_i(\mu(P))=\bar{k}_i(P).
$$
Let $\tau$ denote involution $\sigma\mu$. Involution $\tau$
acts on local parameters as follows
$$
 k_i(\tau(P))=-\bar{k}_i(P).
$$

The following theorem holds:

\begin{theorem}
Let $Q_i$ be fixed points of the antiholomorphic involution
$\tau$. Suppose that on $\Gamma$ there exists a meromorphic 1-form
$\Omega$ with the following set of divisors of zeros and  poles
$$
(\Omega)_0=\gamma+\tau\gamma+P_1+P_2, \eqno{(4)}
$$
$$
 (\Omega)_{\infty}=Q_1+Q_2+Q_3+r+R+\tau R, \eqno{(5)}
$$
where $R=r_1+\dots+r_l$.
Then, function $\varphi^i$ satisfies the equations
$$
 \varphi^1\bar{\varphi}^1A_1+
 \varphi^2\bar{\varphi}^2A_2+\varphi^3\bar{\varphi}^3A_3+
 |d|^2\rm{Res}_{r}\Omega=0,
$$
$$
 \varphi^1\bar{\varphi}^1_xA_1+
 \varphi^2\bar{\varphi}^2_xA_2+
 \varphi^3\bar{\varphi}^3_xA_3=0,
$$
$$
 \varphi^1\bar{\varphi}^1_yA_1+
 \varphi^2\bar{\varphi}^2_yA_2+
 \varphi^3\bar{\varphi}^3_yA_3=0,
$$
$$
 \varphi^1_x\bar{\varphi}^1_yA_1+
 \varphi^2_x\bar{\varphi}^2_yA_2+
 \varphi^3_x\bar{\varphi}^3_yA_3=0,
$$
$$
 \varphi^1_x\bar{\varphi}^1_xA_1+
 \varphi^2_x\bar{\varphi}^2_xA_2+
 \varphi^3_x\bar{\varphi}^3_xA_3-|f_1|^2c_1=0,\eqno{(6)}
$$
$$
 \varphi^1_y\bar{\varphi}^1_yA_1+
 \varphi^2_y\bar{\varphi}^2_yA_2+
 \varphi^3_y\bar{\varphi}^3_yA_3-|f_2|^2c_2=0,\eqno{(7)}
$$
with $A_k=\frac{\rm{Res}_{Q_k}\Omega}{|\alpha_k|^2},\ k=1,2,3$,
$c_1, c_2$ being the coefficients of form $\Omega$ in the
neigborhood of points $P_1$ and $P_2$:
$$
 \Omega=(c_1w_1+i aw_1^2+\dots)dw_1,\ w_1=1/k_1,
$$
$$
 \Omega=(c_2w_2+i bw_2^2+\dots)dw_2,\ w_2=1/k_2.
$$
\end{theorem}

\begin{corollary}
If $\rm{Res}_{Q_i}\Omega>0,\ \rm{Res}_{r}\Omega<0$, then for
$$
 \alpha_i=\sqrt{\rm{Res}_{Q_i}\Omega},\
 d=\sqrt{\frac{-1}{\rm{Res}_{r}\Omega}},
$$
the following equality holds:
$$
   <\varphi,\varphi>=1,\ <\varphi,\varphi_x>=<\varphi,\varphi_y>=<\varphi_x,\varphi_y>=0,
$$
i.e. the mapping ${\cal H}\circ r:{\mathbb R}^2\rightarrow {\mathbb
C}P^2$
is Lagrangian, with the induced metric on $\Sigma$ having a diagonal form
$$
 ds^2=|f_1|^2|c_1|dx^2+|f_2|^2|c_2|dy^2.
$$
\end{corollary}
{\it Proof of theorem 1}. Consider the 1-form
$\Omega_1=\psi(P)\overline{\psi(\tau P)}\Omega$. By virtue of the
definition of involution $\tau$, function $\overline{\psi(\tau
P)}$ has the following form in the neighborhoods of points $P_1$
and $P_2$
$$
 \overline{\psi(\tau P)}=e^{-k_1x}
 \left(\bar{f}_1(x,y)-\frac{\bar{g}_1(x,y)}{k_1}+\frac{\bar{h}_1(x,y)}{k_1^2}+\dots\right),
$$
$$
 \overline{\psi(\tau P)}=e^{-k_2y}\left(\bar{f}_2(x,y)-
 \frac{\bar{g}_2(x,y)}{k_2}+\frac{\bar{h}_2(x,y)}{k_2^2}+\dots\right).
$$
Consequently, form $\Omega_1$ has no essential singularities in
points $P_1$ and $P_2$. The simple poles $\gamma+\tau\gamma$ of
function $\psi(P)\overline{\psi(\tau P)}$ reduce with the zeros in
these points of form $\Omega$. The zeros $R+\tau R$ of function
$\psi(P)\overline{\psi(\tau P)}$ reduce with the simple poles of
form $\Omega$. Thus, form $\Omega_1$ has only simple poles in
points $Q_1,Q_2,Q_3$ and $r$ with equal residues corresponding to
$$
 \psi(Q_1)\overline{\psi(Q_1)}{\rm Res}_{Q_1}\Omega=
 \varphi^1\bar{\varphi}^1A_1,\
 \varphi^2\bar{\varphi}^2A_2,\ \varphi^3\bar{\varphi}^3A_3,\
 |d|^2\rm{Res}_{r}\Omega.
$$
Consequently, the sum of these residues is equal to zero, and this proves the first equality.

Form $\psi(P)\overline{\psi(\tau P)_x}\Omega$ also has no real singularities in points $P_1$ and $P_2$. This form has only simple poles in points $Q_1,Q_2$ and $Q_3$ with residues equal to
$$
 \varphi^1\bar{\varphi}^1_xA_1,\
 \varphi^2\bar{\varphi}^2_xA_2,\
 \varphi^3\bar{\varphi}^3_xA_3.
$$
The second equality is proven. The proof of the following two equalities is analogical.
For this, consider the forms
$$
 \psi(P)\overline{\psi(\tau P)_y}\Omega,\
 \psi(P)_x\overline{\psi(\tau P)_y}\Omega,
$$
that also have only simple poles in points
$Q_1,Q_2$ and $Q_3$.

In order to prove the last two equalities, consider forms
$$
 \psi(P)_x\overline{\psi(\tau P)_x}\Omega,
\
 \psi(P)_y\overline{\psi(\tau P)_y}\Omega.
$$
 These forms have simple poles in points $Q_1,Q_2$, $Q_3, P_1$ and $Q_1,Q_2$, $Q_3, P_2$
with the residues
$$
 \varphi^1_x\bar{\varphi}^1_xA_1,\
 \varphi^2_x\bar{\varphi}^2_xA_2,\
 \varphi^3_x\bar{\varphi}^3_xA_3,\ -|f_1|^2c_1
$$
and
$$
 \varphi^1_y\bar{\varphi}^1_yA_1,\
 \varphi^2_y\bar{\varphi}^2_yA_2,\
 \varphi^3_y\bar{\varphi}^3_yA_3,\ -|f_2|^2c_2.
$$
 Theorem 1 is proven.

{\bf Remark.} {\it Actually, in corollary 1 it is sufficient to demand that only the following inequalities be fulfilled
$$
 \frac{{\rm Res}_{Q_2}\Omega}{{\rm Res}_{Q_1}\Omega}>0,\
 \frac{{\rm Res}_{Q_3}\Omega}{{\rm Res}_{Q_1}\Omega}>0,\
 \frac{{\rm Res}_{r}\Omega}{{\rm Res}_{Q_1}\Omega}<0.
$$
The entire construction easily applies to this case.}

\subsection{Hamiltonian-minimal Lagrangian immersions }
In this section, we find spectral data such that the mapping constructed in the preceding section is conformal and Hamiltonian-minimal.

Consider function
$$
 F(x,y,P)=\partial_x^2\psi+\partial_y^2\psi+A(x,y)\partial_x\psi+
 B(x,y)\partial_y\psi+C(x,y)\psi.
$$
Chose functions $A(x,y),B(x,y)$ and $C(x,y)$ such that
$$
 F(x,y,Q_i)=0,\ i=1,2,3.
$$
If mapping ${\cal H}\circ\varphi$ is conformal, then by lemma 2, the Lagrangian angle is expressed by the functions $A(x,y)$ and $B(x,y)$.
Thus, we need to find spectral data such that the metric on surface $\Sigma$ has a conformal form, and the coefficients $A,B$ are constants.

\begin{lemma} The following equality holds:
$$
 A(x,y)=\frac{1}{c_1|f_1|^2}(-c_1\bar{f}_1(g_1+2{f_1}_x)+f_1(-ia\bar{f}_1+c_1(\bar{g}_1+\bar{f_1}_x))+
 c_2f_2\bar{f_2}_x),
$$
$$
 B(x,y)=\frac{1}{c_2|f_2|^2}(-c_2\bar{f}_2(g_2+2{f_2}_y)+f_2(-ib\bar{f}_2+c_2(\bar{g}_2+\bar{f_2}_y))+
 c_1f_1\bar{f_1}_y),
$$
$$
 C(x,y)=-\frac{1}{|d|^2{\rm Res}_r\Omega}(c_1|f_1|^2+c_2|f_2|^2).
$$
\end{lemma}
{\it Proof of lemma 4.} Consider three forms
$$
 \omega_1=F(P)\overline{\psi(\tau P)}\Omega,
$$
$$
 \omega_2=F(P)\overline{\psi(\tau P)_x}\Omega,
$$
$$
 \omega_3=F(P)\overline{\psi(\tau P)_y}\Omega.
$$
By the construction of involution $\tau$, forms $\omega_1$, $\omega_2$ and $\omega_3$
have no essential singularities in points $P_1$ and $P_2$.

Form $\omega_1$ has only simple poles in points $P_1, P_2$ and
$r$ with the sum of residues
$$
 c_1|f_1|^2+c_2|f_2|^2+C(x,y) |d|^2{\rm Res}_r\Omega =0.
$$
Form $\omega_2$ has only a pole of second order in point $P_1$ and
a simple pole in point $P_2$, whose sum of residues equals
$$
 -c_1\bar{f}_1(g_1+2{f_1}_x)+f_1((-ia-c_1A(x,y))\bar{f_1}+c_1(\bar{g}_1+\bar{f_1}_x))+
 c_2f_2\bar{f_2}_x=0.
$$
Analogically, form $\omega_3$ has only a pole of second order in point $P_2$ and a simple pole in point
$P_1$ with a sum of residues
$$
 -c_2\bar{f}_2(g_2+2{f_2}_y)+f_2((-ib-c_2B(x,y)\bar{f_2}+c_2(\bar{g}_2+\bar{f_2}_y))+
 c_1f_1\bar{f_1}_y=0.
$$
This yields formulas for $A(x,y),B(x,y),C(x,y)$, indicated in
lemma 4. Thus, lemma 4 is proven.

The following lemma holds.

\begin{lemma}
Suppose that on surface $\Gamma$ there exists a
meromorphic form $\omega$ with the following divisors of zeros and poles
$$
 (\omega)_0=\gamma+\sigma\gamma,\ (\omega)_{\infty}=P_1+P_2+R+\sigma
 R,\eqno{(8)}
$$
where ${\rm Res}_{P_1}\omega+{\rm Res}_{P_2}\omega=0$.
Then,
$$
 f_1^2=f_2^2,
$$
i.e. the induced metric on surface $\Sigma$ has the form
$$
 ds^2=|f_1|^2(|c_1|dx^2+|c_2|dy^2).
$$
\end{lemma}
{\it Proof of lemma 5.} Consider the form $\omega_1=\psi(P)\psi(\sigma P)\omega.$ This form has only simple poles in points $P_1$ and $P_2$, consequently, the sum of residues in these points equals zero
$$
 f_1^2{\rm Res}_{P_1}\omega+f_2^2{\rm Res}_{P_2}\omega=0.
$$
Lemma 5 is proven.

{\bf Remark.} {\it Condition ${\rm Res}_{P_1}\omega+{\rm
Res}_{P_2}\omega=0$ is fulfilled, for instance, if $l=0$ or when
$\sigma\omega=-\omega$.}

\begin{lemma}
Suppose that
$$
 \mu(\gamma)=\gamma, \ \mu(R)=R,\ d\in{\mathbb R}.
$$
are real.

Then,
$$
 \psi(x,y,\mu(P))=\overline{\psi(x,y,P)}.
$$
\end{lemma}
For the proof of this standard lemma it is sufficient to note that
function $\overline{\psi(x,y,\mu(P))}$ fulfills conditions 1)--3)
in the definition of the Beker-Akhiezer function. Consequently the
functions $\overline{\psi(x,y,\mu(P))}$ and $\psi(x,y,P)$ are
coin\-cide.

Below we suppose the conditions of lemma 6 to be fulfilled. In particular, this means that
$$
 \tau(\gamma)=\sigma(\gamma), \ \tau(R)=\sigma(R), \eqno{(9)}
$$
and also, that functions $f_i,g_i$, participating in decomposition
$\psi$ in the neigh\-borhood of points $P_1$ and $P_2$ are real.

{\bf Remark.} {\it As follows from the existence of the
meromorphic forms $\omega$ and $\Omega$ on surface $\Gamma$, there
is a meromorphic function $\lambda=\frac{\omega}{\Omega}$ with the
following divisors of zeros and poles
$$
 (\lambda)_0=Q_1+Q_2+Q_3+r,\ (\lambda)_{\infty}=2P_1+2P_2
$$}
on $\Gamma$.

To begin with, suppose that in the decomposition of form $\Omega$,
in neighborhoods of points $P_1$ the coefficient $c_1$ is equal to
1. By virtue of forms $\Omega$ and $\overline{\Omega(\tau P)}$
possessing the same set of zeros and poles, forms $\Omega$ and
$\overline{\Omega(\tau P)}$ are proportional. From $\tau
w_1=-\overline{w}_1$ and from the fact that in the neighborhood of
point $P_1$ decompositions
$$
 \Omega(P)=(w_1+i aw_1^2+\dots)dw_1,
$$
$$
 \overline{\Omega(\tau P)}=(w_1+i \bar{a}w_1^2+\dots)dw_1,
$$
hold,
we get
$$
 \Omega(P)=\overline{\Omega(\tau P)}.
$$
In the neighborhood of point $P_2$ we have
$$
 \Omega(P)=(c_2w_2+i bw_2^2\dots)dw_2=\overline{\Omega(\tau P)}=(\bar{c}_2w_2+i \bar{b}w_2^2\dots)dw_2,
$$
consequently, $c_2,a,b$ are real constants.

Suppose that the inequalities of residues in corollary 1 be fulfilled. From (6) and (7) follows that $c_1$ and $c_2$ have the same sign, and by virtue of assumption $c_1=1$, $c_2$ is a positive constant. Thus, replacing
$w_2\rightarrow \frac{w_2}{\sqrt{c_2}}$, consider that
$$
 c_1=c_2=1.
$$
Thus, the metric on surface $\Sigma$ has a conformal form
$$
 ds^2=f_1^2(dx^2+dy^2).
$$

The main theorem holds.

\begin{theorem}
Suppose that the conditions of theorem 1 and lemmas 5 and 6 are
fulfilled. Then, the Lagrangian angle looks as follows
$$
 \beta= a x+ b y+ c,
$$
where $c$ is some real constant, i.e. surface $\Sigma$ is Hamiltonian-minimal.
If also $\sigma Q_i=Q_i$, then surface $\Sigma$ is minimal.
\end{theorem}
{\it Proof of theorem 2.}
Since
$$
 f_1=\bar{f_1},\ f_2=\bar{f_2},\ f_1^2=f_2^2, \ g_1=\bar{g}_1,\ g_2=\bar{g}_2,\ c_1=c_2=1
$$
by lemma 4, coefficients $A$ and $B$ have the form
$$
 A=-i a,\ B=-i b,
$$
consequently,  the Lagrangian angle has the form indicated in theorem 2, consequently, surface $\Sigma$ is Hamiltonian-minimal.

What's left is to prove the last point of the theorem.
Consider two forms
$\Omega(P)$ and $\Omega(\sigma P)$. According to our assumptions, these forms have the same set of zeros and poles.
From their decompositions in the neighborhood of points $P_1$ and $P_2$
$$
 \Omega(P)=(w_1+i aw_1^2+\dots)dw_1,
$$
$$
 \Omega(\sigma P)=(w_1-i aw_1^2+\dots)dw_1,
$$
$$
 \Omega(P)=(w_2+i bw_2^2+\dots)dw_2,
$$
$$
 \Omega(\sigma P)=(w_2-i bw_2^2+\dots)dw_2
$$
follows that forms $\Omega(P)$ and $\Omega(\sigma P)$ coincide and $a=b=0$, i.e. surface
$\Sigma$ is minimal.
Thus, theorem 2 is proven.

\section{Examples}
\subsection{Spectral curve of genus 0: Clifford tori}
Suppose that $\Gamma={\mathbb C}P^1, P_1=\infty, P_2=0$. The Riemannian surface $\Gamma$ allows the holomorphic involution
$$
 \sigma: \Gamma\rightarrow\Gamma,\ \sigma(z)=-z
$$
and the antiholomorphic involutions
$$
 \mu: \Gamma\rightarrow\Gamma,\ \mu(z)=\bar{z},
$$
$$
 \tau: \Gamma\rightarrow\Gamma,\ \tau(z)=-\bar{z}.
$$
Let $l=0$, then the Baker-Akhiezer function has the form
$$
 \psi(x,y,z)=e^{xz+\frac{y}{z}}f(x,y).
$$
Let points $Q_1, Q_2, Q_3$ have the coordinates $iq_1, iq_2, iq_3\in
i{\mathbb R}$, and point $r$  --- a coordinate $i\tilde{r}\in
i{\mathbb R}$. In this case $\tau(Q_j)=Q_j,\ \tau(r)=r.$ From the normalization condition follows $\psi(x,y,r)=d.$ This leads us to
$$
 \psi(x,y,z)=de^{xz+\frac{y}{z}}e^{-i\tilde{r}x+\frac{iy}{\tilde{r}}}.
$$
The form $\Omega$ looks as follows:
$$
 \Omega=\frac{szdz}{(z-iq_1)(z-iq_2)(z-iq_3)(z-i\tilde{r})}.
$$
The decomposition coefficients $c_1$, $c_2$ of form $\Omega$ in the neighborhood of points
$P_1$, $P_2$ equal
$$
 c_1=-s,\ c_2=\frac{s}{q_1q_2q_3\tilde{r}}.
$$
Consequently,
$$
 s=-1,\ \tilde{r}=-\frac{1}{q_1q_2q_3}.
$$
The vector-function $\varphi$ has the following components
$$
 \varphi^1=\alpha_1de^{\frac{-i(q_1-\tilde{r})(q_1\tilde{r}x+y)}{q_1\tilde{r}}},\
 \varphi^2=\alpha_2de^{\frac{-i(q_2-\tilde{r})(q_2\tilde{r}x+y)}{q_2\tilde{r}}},
\
 \varphi^3=\alpha_3de^{\frac{-i(q_3-\tilde{r})(q_3\tilde{r}x+y)}{q_3\tilde{r}}},
$$
where
$$
 \alpha_1=\sqrt{\frac{-q_1}{(q_2-q_1)(q_1-q_3)(q_1-\tilde{r})}},\
 \alpha_2=\sqrt{\frac{-q_2}{(q_1-q_2)(q_2-q_3)(q_2-\tilde{r})}},\
$$
$$
 \alpha_3=\sqrt{\frac{-q_3}{(q_1-q_3)(q_3-q_2)(q_3-\tilde{r})}},\
 d=\sqrt{\frac{(\tilde{r}-q_1)(\tilde{r}-q_2)(\tilde{r}-q_3)}{-\tilde{r}}}.
$$
For example, for $q_1=2, q_2=-2,q_3=1/2$ we have
$$
 \varphi^1=\frac{3}{\sqrt{14}}e^{\frac{1}{2}i(2x+y)},\
 \varphi^2=\frac{1}{\sqrt{6}}e^{\frac{3}{2}i(-2x+y)},
\
 \varphi^3=\frac{2}{\sqrt{21}}e^{-3i(\frac{x}{4}+y)}.
$$
The induced metric looks as follows:
$$
 ds^2=\frac{7}{4}(dx^2+dy^2),
$$
with functions $\varphi^j$ satisfying equation
$$
 \Delta \varphi^j+i\frac{11}{4}\partial_x
 \varphi^j+i\partial_y\varphi^j+\frac{9}{2}\varphi^j=0,
$$
the Lagrangian angle has the form
$$
 \beta=-\left(\frac{11}{4}x+y-\pi\right).
$$
\subsection{Spectral curves of genus $g>0$: {\bf HML}-tori}
Let $\Gamma_0$ be a hyperelliptic surface of genus $g$
given by equation
$$
 y^2=P(x),
$$
where $P(x)$ is a polynom of degree $2g+2$ with real coefficients without multiple roots. Let $f$ denote a meromorphic function on surface $\Gamma_0$
$$
 f=\frac{x-\beta}{x-\alpha},
$$
where $\alpha, \beta$  are some real numbers such that $P(\alpha)\ne 0, P(\beta)\ne 0$. Function $f$ has two simple zeros in points
$$
 \tilde{P_1}=(\beta,\sqrt{P(\beta )}),\
 \tilde{P_2}=(\beta,-\sqrt{P(\beta)})
$$
and two simple poles in points
$$
 \tilde{r}=(\alpha,\sqrt{P(\alpha)}),\
 \tilde{Q}_1=(\alpha,-\sqrt{P(\alpha)}).
$$
Let $\Gamma$ denote the Riemannian surface of function $\sqrt{f}$.
The affine part of surface $\Gamma$ is given in ${\mathbb C}^3$ with
coordinates $x,y,z$ by the system of equations
$$
 y^2=P(x),
$$
$$
 z^2=\frac{x-\beta}{x-\alpha}.
$$
Surface $\Gamma$ allows a holomorphic involution
$$
 \sigma:\Gamma\rightarrow \Gamma,\ \sigma(x,y,z)=(x,y,-z)
$$
with four fixed points $P_1,P_2,Q_1,r$ --- the inverse images of points $\tilde{P_1},\tilde{P_2},\tilde{Q}_1,\tilde{r}$ for the
projection $\Gamma\rightarrow\Gamma_0$.

Surface $\Gamma$ also allows antiholomorphic involutions
$$
 \mu:\Gamma\rightarrow \Gamma,\
 \mu(x,y,z)=(\bar{x},\bar{y},\bar{z}),
$$
$$
 \tau:\Gamma\rightarrow \Gamma,\
 \tau(x,y,z)=(\bar{x},\bar{y},-\bar{z}).
$$
Choose $\alpha$ and $\beta$ such that $P(\alpha)>0, P(\beta)> 0$.
Then, points $P_1,P_2,Q_1,r$ are fixed with respect to involution
$\tau$.

Consider the meromorphic function $\lambda$ on surface $\Gamma$
$$
 \lambda=z^2-c z,
$$
where $c$ is some imaginary constant. This function has two multiple poles in points $P_1,P_2$, two simple zeros in points $Q_1, r$ and two zeros in points
$$
 Q_2=\left(\frac{\beta-c^2\alpha}{1-c^2},\sqrt{P\left(\frac{\beta-c^2\alpha}{1-c^2}\right)},c\right),
$$
$$
 Q_3=\left(\frac{\beta-c^2\alpha}{1-c^2},-\sqrt{P\left(\frac{\beta-c^2\alpha}{1-c^2}\right)},c\right).
$$
Choose $c$ to be purely imaginary such that
$$
 P\left(\frac{\beta-c^2\alpha}{1-c^2}\right)>0.
$$
Then, points $Q_2$ and $Q_3$ are invariant with respect to $\tau$, but not invariant with respect to $\sigma$.

Thus, if form $\omega$ has a divisor of zeros and poles of form
 (8) with property $\sigma \omega=-\omega$ with
real divisors $\gamma$ and $R$ (see (9)), then form $\Omega=\frac{\omega}{\lambda}$ has the required set of zeros and poles (4) and (5).

If needed, divide $\Omega$ by a suitable constant such that the decomposition $\Omega$ in the neighborhood of $P_1$ has the form
$$
 \Omega=(w_1+\dots)dw_1.
$$
Forms $\overline{\Omega(\tau P)}$ and $\Omega(P)$ have the same
set of zeros and poles. Consequently, by virtue of decomposition
$\Omega$ in the vicinity of $P_1$ and by virtue of $\tau
(w_1)=-\bar{w}_1$ get $\overline{\Omega(\tau P)}=\Omega(P)$. Chose
a local parameter $v_i$ in the neighborhood of point $Q_i$ such
that $\tau(v_i)=-\bar{v}_i$, consequently, the decomposition
coefficient $q_i$
$$
 \Omega=\left(\frac{q_i}{v_i}+\dots\right)dv_i,
$$
by virtue of equality $\overline{\Omega(\tau P)}=\Omega(P)$, is real, i.e.
$$
 {\rm Res}_{Q_i}\Omega\in{\mathbb R},\ i=1,2,3.
$$
Analogically,
$$
 {\rm Res}_{r}\Omega\in{\mathbb R}.
$$

Now chose the construction's parameters such that
$$
 {\rm Res}_{Q_i}\Omega>0,\ {\rm Res}_{r}\Omega<0
$$
and obtain spectral data satisfying all requirements of theorem 2.

\subsection{Spectral curves of genus $g>0$: minimal Lagrangian tori}
As in the preceding example, let $\Gamma_0$ denote a hyperelliptic
curve of genus $g$ given by equation
$$
 y^2=P(x).
$$
Let $\Gamma$ denote a Riemannian surface of function $\sqrt{g}$, where
$$
 g=\frac{(x-\beta_1)(x-\beta_2)}{(x-\alpha_1)(x-\alpha_2)},\
 \beta_1,\beta_2\in{\mathbb R},
$$
where $\alpha_1,\alpha_2$ are real roots of the polynomial $P(x)$.
The affine part of surface $\Gamma$ is given in ${\mathbb C}^3$ with
coordinates $x,y,z$ by the system of equations
$$
 y^2=P(x),
$$
$$
 z^2=\frac{(x-\beta_1)(x-\beta_2)}{(x-\alpha_1)(x-\alpha_2)}.
$$
Let $P_1,P_2\in\Gamma$ denote inverse images of points $(\alpha_1,0),
(\alpha_2,0)\in\Gamma_0$ for a natural projection
$\Gamma\rightarrow \Gamma_0$, let $Q_1,Q_2,Q_3,r\in\Gamma$
denote points
$$
 Q_1=(\beta_1,\sqrt{P(\beta_1)},0),\ Q_2=(\beta_1,-\sqrt{P(\beta_1)},0),
$$
$$
 Q_3=(\beta_2,\sqrt{P(\beta_2)},0),\
 r=(\beta_2,-\sqrt{P(\beta_2)},0).
$$
All these points are invariant under a holomorphic involution
$$
 \sigma:\Gamma\rightarrow \Gamma,\ \sigma(x,y,z)=(x,y,-z).
$$
Surface $\Gamma$ allows antiholomorphic involutions
$$
 \mu:\Gamma\rightarrow \Gamma,\
 \mu(x,y,z)=(\bar{x},\bar{y},\bar{z}),
$$
$$
 \tau:\Gamma\rightarrow \Gamma,\
 \tau(x,y,z)=(\bar{x},\bar{y},-\bar{z}).
$$
Chose $\beta_1$ and $\beta_2$ such that $P(\beta_1)>0, P(\beta_2)> 0$.
Then, points $P_1,P_2,Q_1,$ $Q_2$, $Q_3,r$ are fixed for involution $\tau$.
Function $\lambda=z$ has the following set of zeros and poles
$$
 (\lambda)_0=Q_1+Q_2+Q_3+r,\ (\lambda)_{\infty}=2P_1+2P_2.
$$
Thus, if form $\omega$ has simple poles in points
$P_1$ and $P_2$ and $\sigma\omega=-\omega$, then form $\Omega=\frac{\omega}{\lambda}$
has a set of zeros and poles satisfying the condition of theorem 2. Since $\sigma(Q_i)=Q_i$,
according to theorem 2 surface $\Sigma$ is minimal.

$$
$$

\noindent Sobolev Institute of Mathematics,

\noindent Novosibirsk State University,

\noindent {\it E-mail address:} mironov@math.nsc.ru

\end{document}